\input amstex
\documentstyle{amsppt}\nologo\footline={}

\def\Gr{\mathop{\text{\rm Gr}}}
\def\Lin{\mathop{\text{\rm Lin}}}
\def\T{\mathop{\text{\rm T}}}
\def\tr{\mathop{\text{\rm tr}}}
\def\GL{\mathop{\text{\rm GL}}}
\def\Re{\mathop{\text{\rm Re}}}
\def\ricci{\mathop{\text{\rm ricci}}}
\def\G{\mathop{\text{\rm G}}}
\def\Span{\mathop{\text{\rm Span}}}
\def\B{\mathop{\text{\rm B}}}

\hsize450pt

\topmatter\title Differential geometry of grassmannians and
Pl\"ucker map\endtitle\author Sasha Anan$'$in and Carlos
H.~Grossi\endauthor\thanks First author
partially supported by the Institut des Hautes \'Etudes Scientifiques
(IH\'ES).\endthanks\thanks Second author supported by the
Max-Planck-Gesellschaft.\endthanks\address Departamento de
Matem\'atica, IMECC, Universidade Estadual de
Campinas,\newline13083-970--Campinas--SP, Brasil\endaddress\email
Ananin$_-$Sasha\@yahoo.com\endemail\address Max-Planck-Institut f\"ur
Mathematik, Vivatsgasse 7, 53111 Bonn, Germany\endaddress\email
grossi$_-$ferreira\@yahoo.com\endemail\subjclass 53A20 (53A35,
51M10)\endsubjclass\abstract Using the Pl\"ucker map between
grassmannians, we study basic aspects of classic grassmannian
geometries. For `hyperbolic' grassmannian geometries, we prove some
facts (for instance, that the Pl\"ucker map is a minimal isometric
embedding) that were previously known in the `elliptic'
case.\endabstract\endtopmatter\document

\centerline{\bf1.~Introduction}

\medskip

We study the nondegenerate piece $\Gr^0(k,V)$ of the grassmannian
$\Gr(k,V)$ of $k$-dimensional subspaces in an $\Bbb R$- or
$\Bbb C$-vector space $V$ equipped with an hermitian form. This paper
links the (pseudo-)riemannian geometry of $\Gr^0(k,V)$ to structures
discussed in [AGr] and [AGoG]. It is merely intended to illustrate how
do the methods from the mentioned papers work in the differential
geometry of grassmannians. Many of the results presented here are known
in particular cases.\footnote{If the hermitian form on $V$ is definite,
the classic geometry is sort of elliptic. Most of the known facts deal
with this case. The hermitian algebra of the indefinite form requires
additional effort thus making it nontrivial the case of `hyperbolic'
classic geometries.}
We believe that our treatment provides additional clarity even in those
cases.

It follows a brief description of the results. The Pl\"ucker map is a
minimal isometric embedding. The Gauss equation provides the curvature
tensor in the form of the $(2,1)$-symmetrization of the triple product
exactly as in the projective case [AGr, Subsection 4.4]. $\Gr^0(k,V)$
is shown to be Einstein. Generic geodesics in $\Gr^0(k,V)$ are
described. Also, we illustrate how a grassmannian classic geometry
unexpectedly shows up in relation to convexity in real hyperbolic
space.

It turns out that the hermitian metric actually plays no role in most
of the proofs. The tangent vectors can usually be taken as {\it
footless\/} or as {\it observed\/} from different points. Therefore,
many concepts, for instance, those of isometric or minimal embeddings
and of the Gauss equation, may be restated in the terms of the {\it
product\/} (see [AGr, Subsection 1.1] or [AGoG, Sections 2, 3] for the
definitions). This must be fruitful since the product embodies
different (pseudo-)riemannian concepts in a single simple structure. In
the spirit of [AGoG], it would be nice to understand what remains from
these concepts after arriving at the absolute.

To prevent a possible scepticism of the reader, we have to say that the
pseudo-riemannian metrics play a fundamental role in the study of the
riemannian classical geometries: basic geometrical objects almost
never form riemannian spaces. To illustrate this remark, the beautiful
article [GuK] is to be mentioned, where the authors work in an ambient
that in fact falls into our settings.

The differential geometry of grassmannians is a rather vast field (see,
for instance, the survey [BoN]). We believe that it is reasonable to
redemonstrate known facts in the area by using the language of [AGr]
and [AGoG]. Of course, we recognize that such a project involves a huge
amount of work, but is probably worth the candle: besides giving each
fact an appropriate generality, it would provide a better understanding
of particular problems in classic geometries.

\bigskip

\centerline{\bf2.~Pl\"ucker-and-play}

\medskip

We remind some notation and convention from [AGoG, Section 2]. Let $V$
be an $n$-dimensional $\Bbb K$-vector space equipped with a
nondegenerate hermitian form $\langle-,-\rangle$, where $\Bbb K=\Bbb R$
or $\Bbb K=\Bbb C$. Take and fix a $\Bbb K$-vector space $P$ such that
$\dim_\Bbb KP=k$. Denote by
$M:=\big\{p\in\Lin_\Bbb K(P,V)\mid\ker p=0\big\}$ the open subset of
all monomorphisms in $\Lin_\Bbb K(P,V)$. The group $\GL_\Bbb KP$ acts
from the right on $\Lin_\Bbb K(P,V)$ and on $M$. The~grassmannian
is the quotient space $\pi:M\to\Gr_\Bbb K(k,V):=M/\GL_\Bbb KP$. We do
not distinguish between the notation of points in $\Gr_\Bbb K(k,V)$ and
of their representatives in $M$. We~frequently write $p$ in place of
the image $pP$ and $p^\perp$, in place of the orthogonal $(pP)^\perp$.
The space $\Gr_\Bbb K^0(k,V)\subset\Gr_\Bbb K(k,V)$ is formed by the
nondegenerate subspaces. The tangent space $\T_pM$ is commonly
identified with $\Lin_\Bbb K(p,V)$. For~$p\in\Gr_\Bbb K^0(k,V)$, we
identify
$\T_p\Gr_\Bbb K^0(k,V)=\Lin_\Bbb K(p,p^\perp)\subset\Lin_\Bbb K(V,V)$,
where the inclusion is provided by $V=p\oplus p^\perp$.

Our purpose is to study the {\it$m$-Pl\"ucker embedding\/}
$$\textstyle E^m:\Gr_\Bbb K(k,V)\to\Gr_\Bbb K\big({k\choose
m},\bigwedge^mV\big),\qquad p\mapsto\bigwedge^mp,$$
where the vector space $\bigwedge^mV$ is equipped with the hermitian
form given by the rule
$$\langle v_1\wedge\dots\wedge v_m,w_1\wedge\dots\wedge
w_m\rangle:=\det\langle v_i,w_j\rangle.$$

Let $p\in M$. It is not difficult to see that the differential of the
map $M\to\Lin_\Bbb K\big(\bigwedge^mP,\bigwedge^mV\big)$ at $p$ sends
the tangent vector $\overline t\in\T_pM=\Lin_\Bbb K(p,V)$ to
$E^m\overline t\in\Lin_\Bbb K\big(\bigwedge^mp,\bigwedge^mV\big)$
defined by the rule
$E^m\overline t:p_1\wedge\dots\wedge
p_m\mapsto\sum_{i=1}^mp_1\wedge\dots\wedge\overline
tp_i\wedge\dots\wedge p_m$
for all $p_1,\dots,p_m\in p$. Therefore, we can describe the
differential of $E^m$ at $p$ as
$$\textstyle E^m:\Lin_\Bbb K(p,V/p)\to\Lin_\Bbb
K\big(\bigwedge^mp,\bigwedge^mV/\bigwedge^mp\big),$$
$$\textstyle E^mt:p_1\wedge\dots\wedge
p_m\mapsto{\displaystyle\sum_{i=1}^m}p_1\wedge\dots\wedge\overline
tp_i\wedge\dots\wedge p_m+\bigwedge^mp$$
for all $t:p\to V/p$ and $p_1,\dots,p_m\in p$, where
$\overline t:p\to V$ is an arbitrary lift of $t$.

Given $p\in\Gr_\Bbb K^0(k,V)$, we have the orthogonal decomposition
$$\textstyle\bigwedge^mV={\displaystyle\bigoplus_{i=0}^m}\bigwedge^i
p^\perp\wedge \bigwedge^{m-i}p.\leqno{\bold{(2.1)}}$$
In particular, taking $p\in\Gr_\Bbb K^0(k,V)$ and
$t\in\T_p\Gr_\Bbb K^0(k,V)=\Lin_\Bbb K(p,p^\perp)$, we obtain
$$E^mt:p_1\wedge\dots\wedge p_m\mapsto\sum_{i=1}^mp_1\wedge\dots\wedge
tp_i\wedge\dots\wedge p_m\leqno{\bold{(2.2)}}$$
for all $p_1,\dots,p_m\in p$. Note that (2.2) makes sense for an
arbitrary $t:V\to V$.

Define the linear map
$B(t_1,t_2):\bigwedge^mV\to\bigwedge^mV$ by the rule
$$B(t_1,t_2)(v_1\wedge\dots\wedge v_m):=\sum_{i\ne
j}v_1\wedge\dots\wedge t_1v_i\wedge\dots\wedge t_2v_j\wedge\dots\wedge
v_m$$
for all $v_1,\dots,v_m\in V$, where $t_1,t_2:V\to V$. (In the above
sum, $t_2v_j$ appears before $t_1v_i$ if $i>j$.)

\medskip

{\bf2.3.~Lemma.} {\sl Let\/ $p\in\Gr_\Bbb K(k,V)$ and let\/
$t,t_1,t_2:V\to V$. Then
$$\big\langle E^mt(p\wedge\dots\wedge p_m),q\wedge
v_2\wedge\dots\wedge v_m\big\rangle=\langle p_1\wedge\dots\wedge
p_m,t^*q\wedge v_2\wedge\dots\wedge v_m\rangle,$$
$$\big\langle B(t_1,t_2)(p_1\wedge\dots\wedge p_m),q_1\wedge q_2\wedge
v_3\wedge\dots\wedge v_m\big\rangle=$$
$$=\langle p_1\wedge\dots\wedge p_m,t_1^*q_1\wedge t_2^*q_2\wedge
v_3\wedge\dots\wedge v_m\rangle+\langle p_1\wedge\dots\wedge
p_m,t_2^*q_1\wedge t_1^*q_2\wedge v_3\wedge\dots\wedge v_m\rangle$$
for all\/ $q,q_1,q_2\in p^\perp$, $p_1,\dots,p_m\in p$, and
$v_2,\dots,v_m\in V$.}

\medskip

{\bf Proof} is based on simple known identities involving determinants
(marked with $\dag$ and left without proof). We have
$$\big\langle E^mt(p_1\wedge\dots\wedge p_m),q\wedge
v_2\wedge\dots\wedge
v_m\big\rangle=\sum_{i=1}^m\det\left(\smallmatrix0&\langle
p_1,v_2\rangle&\cdots&\langle p_1,v_m\rangle\\
\vdots&\vdots&\ddots&\vdots\\0&\langle
p_{i-1},v_2\rangle&\cdots&\langle p_{i-1},v_m\rangle\\\langle
tp_i,q\rangle&\langle tp_i,v_2\rangle&\cdots&\langle
tp_i,v_m\rangle\\0&\langle p_{i+1},v_2\rangle&\cdots&\langle
p_{i+1},v_m\rangle\\\vdots&\vdots&\ddots&\vdots\\0&\langle
p_m,v_2\rangle&\cdots&\langle
p_m,v_m\rangle\endsmallmatrix\right)\overset\dag\to=$$
$$\overset\dag\to=\det\left(\smallmatrix\langle tp_1,q\rangle&\langle
p_1,v_2\rangle&\cdots&\langle
p_1,v_m\rangle\\\vdots&\vdots&\ddots&\vdots\\\langle
tp_m,q\rangle&\langle p_m,v_2\rangle&\cdots&\langle
p_m,v_m\rangle\endsmallmatrix\right)=\det\left(\smallmatrix\langle
p_1,t^*q\rangle&\langle p_1,v_2\rangle&\cdots&\langle
p_1,v_m\rangle\\\vdots&\vdots&\ddots&\vdots\\\langle
p_m,t^*q\rangle&\langle p_m,v_2\rangle&\cdots&\langle
p_m,v_m\rangle\endsmallmatrix\right){=}\langle p_1\wedge\dots\wedge
p_m,t^*q\wedge v_2\wedge\dots\wedge v_m\rangle$$
and
$$\big\langle B(t_1,t_2)(p_1\wedge\dots\wedge p_m),q_1\wedge q_2\wedge
v_3\wedge\dots\wedge v_m\big\rangle=\sum\limits_{i\ne
j}\det\left(\smallmatrix0&0&\langle p_1,v_3\rangle&\cdots&\langle
p_1,v_m\rangle\\\vdots&\vdots&\vdots&\ddots&\vdots\\0&0&\langle
p_{i-1},v_3\rangle&\cdots&\langle p_{i-1},v_m\rangle\\\langle
t_1p_i,q_1\rangle&\langle t_1p_i,q_2\rangle&\langle
t_1p_i,v_3\rangle&\cdots&\langle t_1p_i,v_m\rangle\\0&0&\langle
p_{i+1},v_3\rangle&\cdots&\langle
p_{i+1},v_m\rangle\\\vdots&\vdots&\vdots&\ddots&\vdots\\0&0&\langle
p_{j-1},v_3\rangle&\cdots&\langle p_{j-1},v_m\rangle\\\langle
t_2p_j,q_1\rangle&\langle t_2p_j,q_2\rangle&\langle
t_2p_j,v_3\rangle&\cdots&\langle t_2p_j,v_m\rangle\\0&0&\langle
p_{j+1},v_3\rangle&\cdots&\langle
p_{j+1},v_m\rangle\\\vdots&\vdots&\vdots&\ddots&\vdots\\0&0&\langle
p_m,v_3\rangle&\cdots&\langle
p_m,v_m\rangle\endsmallmatrix\right)\overset\dag\to=$$
$$\overset\dag\to=\det\left(\smallmatrix\langle
t_1p_1,q_1\rangle&\langle t_2p_1,q_2\rangle&\langle
p_1,v_3\rangle&\cdots&\langle
p_1,v_m\rangle\\\vdots&\vdots&\vdots&\ddots&\vdots\\\langle
t_1p_m,q_1\rangle&\langle t_2p_m,q_2\rangle&\langle
p_m,v_3\rangle&\cdots&\langle
p_m,v_m\rangle\endsmallmatrix\right)+\det\left(\smallmatrix\langle
t_2p_1,q_1\rangle&\langle t_1p_1,q_2\rangle&\langle
p_1,v_3\rangle&\cdots&\langle
p_1,v_m\rangle\\\vdots&\vdots&\vdots&\ddots&\vdots\\\langle
t_2p_m,q_1\rangle&\langle t_1p_m,q_2\rangle&\langle
p_m,v_3\rangle&\cdots&\langle p_m,v_m\rangle\endsmallmatrix\right)=$$
$$=\det\left(\smallmatrix\langle p_1,t_1^*q_1\rangle&\langle
p_1,t_2^*q_2\rangle&\langle p_1,v_3\rangle&\cdots&\langle
p_1,v_m\rangle\\\vdots&\vdots&\vdots&\ddots&\vdots\\\langle
p_m,t_1^*q_1\rangle&\langle p_m,t_2^*q_2\rangle&\langle
p_m,v_3\rangle&\cdots&\langle
p_m,v_m\rangle\endsmallmatrix\right)+\det\left(\smallmatrix\langle
p_1,t_2^*q_1\rangle&\langle p_1,t_1^*q_2\rangle&\langle
p_1,v_3\rangle&\cdots&\langle
p_1,v_m\rangle\\\vdots&\vdots&\vdots&\ddots&\vdots\\\langle
p_m,t_2^*q_1\rangle&\langle p_m,t_1^*q_2\rangle&\langle
p_m,v_3\rangle&\cdots&\langle p_m,v_m\rangle\endsmallmatrix\right)=$$
$$=\langle p_1\wedge\dots\wedge p_m,t_1^*q_1\wedge t_2^*q_2\wedge
v_3\wedge\dots\wedge v_m\rangle+\langle p_1\wedge\dots\wedge
p_m,t_2^*q_1\wedge t_1^*q_2\wedge v_3\wedge\dots\wedge v_m\rangle\
_\blacksquare$$

\medskip

Let $t\in\Lin_\Bbb K(p,p^\perp)\subset\Lin_\Bbb K(V,V)$. It follows
from (2.2) and Lemma 2.3 that the only nonvanishing component of
$(E^mt)^*$ related to the decomposition (2.1) has the form
$(E^mt)^*:p^\perp\wedge\bigwedge^{m-1}p\to\bigwedge^mp$,
$$(E^mt)^*:q\wedge p_2\wedge\dots\wedge p_m\mapsto t^*q\wedge
p_2\wedge\dots\wedge p_m,\leqno{\bold{(2.4)}}$$
where $q\in p^\perp$ and $p_2,\dots,p_m\in p$. In other words,
$(E^mt)^*=E^mt^*$. Similar arguments are applicable to $B(t_1,t_2)$
with $t_1,t_2\in\Lin_\Bbb K(p,p^\perp)\subset\Lin_\Bbb K(V,V)$.

\medskip

{\bf2.5.~Proposition {\rm(compare to [BoN, Assertions 1--2])}.} {\it
The\/ $m$-Pl\"ucker embedding provides an hermitian\/ {\rm(}hence,
pseudo-riemannian\/{\rm)} embedding\/
$E^m:\Gr_\Bbb K^0(k,V)\to\Gr_\Bbb K^0\big({k\choose
m},\bigwedge^mV\big)$, assuming the metric on\/ $\Gr_\Bbb K^0(k,V)$
rescaled by the factor\/ ${{k-1}\choose{m-1}}$.}

\medskip

{\bf Proof.} Let $p\in\Gr_\Bbb K^0(k,V)$ and let $t_1,t_2:p\to p^\perp$
be tangent vectors at $p$. By (2.2) and (2.4),
$$(E^mt_1)^*E^mt_2:p_1\wedge\dots\wedge
p_m\mapsto\sum_{i=1}^mp_1\wedge\dots\wedge t_1^*t_2p_i\wedge\dots\wedge
p_m$$
for all $p_1,\dots,p_m\in p$. As is easy to see,
$\tr(E^m\varphi)={{k-1}\choose{m-1}}\tr\varphi$ for every linear map
$\varphi:p\to p$ and the map $E^m\varphi:\bigwedge^mp\to\bigwedge^mp$
defined as in (2.2). Hence,
$$\langle E^mt_1,E^mt_2\rangle=\tr\big((E^mt_1)^*E^mt_2\big)=
\tr\big(E^m(t_1^*t_2)\big)=\textstyle{{k-1}\choose{m-1}}\tr(t_1^*t_2)=
{{k-1}\choose{m-1}}\langle t_1,t_2\rangle\ _\blacksquare$$

\medskip

Given $p\in\Gr_\Bbb K^0(k,V)$, denote by $\pi'[p]$ and $\pi[p]$ the
orthogonal projectors corresponding to the decomposition
$V=p\oplus p^\perp$. For $t\in\Lin_\Bbb K(V,V)$, define the tangent
vector $t_p:=\pi[p]t\pi'[p]$ at $p$.

Let $U\subset M$ be a {\it saturated\/} and {\it nondegenerate\/} open
set. This means that $U\GL_\Bbb KP=U$ and
$\pi U\subset\Gr_\Bbb K^0(k,V)$, where $\pi:M\to\Gr_\Bbb K(k,V)$ stands
for the quotient map. A smooth map $X:U\to\Lin_\Bbb K(V,V)$ is said to
be a {\it lifted field\/} over $U$ if $X(p)_p=X(p)$ and $X(pg)=X(p)$
for all $p\in U$ and $g\in\GL_\Bbb KP$. In~other words, $\pi$ maps $X$
onto a correctly defined smooth tangent field over the open subset
$\pi U\subset\Gr_\Bbb K^0(k,V)$.

For $t\in\Lin_\Bbb K(V,V)$, define
$$\nabla_tX(p):=\Big(\displaystyle\frac
d{d\varepsilon}\Big|_{\varepsilon=0}X\big((1+\varepsilon
t)p\big)\Big)_p.$$
Since $\pi'[pg]=\pi'[p]$ and $\pi[pg]=\pi[p]$ for all $p\in U$ and
$g\in\GL_\Bbb KP$, the field $p\mapsto\nabla_{Y(p)}X$ is lifted for
arbitrary lifted fields $X$ and $Y$ over $U$. Obviously, $\nabla$
enjoys the properties of an affine connection; we~assume
$\Gr_\Bbb K^0(k,V)$ equipped with this {\it intrinsic\/} connection.

\medskip

{\bf2.6.~Proposition.} {\sl The connection induced by the\/
$m$-Pl\"ucker embedding coincides with the intrinsic one and the map
$$B(t_1,t_2):{\T}_p{\Gr}_\Bbb K^0(k,V)\times{\T}_p{\Gr}_\Bbb K^0(k,V)
\to\big(E^m{\T}_p{\Gr}_\Bbb K^0(k,V)\big)^\perp$$
is the second fundamental form of the embedding.}

\medskip

{\bf Proof.} Let $p\in\Gr_\Bbb K^0(k,V)$ and let
$t\in\Lin_\Bbb K(p,p^\perp)\subset\Lin_\Bbb K(V,V)$. First, we need to
establish some auxiliary formulae.

Denote $g(\varepsilon):=1+\varepsilon t$. We have
$\displaystyle\frac
d{d\varepsilon}\Big|_{\varepsilon=0}g(\varepsilon)=t$
and
$\displaystyle\frac d{d\varepsilon}\Big|_{\varepsilon=0}\big(g^{-1}
(\varepsilon)\big)^*=-t^*$
because $g^{-1}(\varepsilon)g(\varepsilon)=1$ for small $\varepsilon$.
The projectors
$$\pi'(\varepsilon):=\pi'\big[{\textstyle{\bigwedge}^m}g(\varepsilon)p
\big],\quad\pi(\varepsilon):=\pi\big[{\textstyle{\bigwedge}^m}g(
\varepsilon)p\big]$$
satisfy
$$\pi'(\varepsilon)\big(g(\varepsilon)p_1\wedge\dots\wedge
g(\varepsilon)p_m\big)=g(\varepsilon)p_1\wedge\dots\wedge
g(\varepsilon)p_m,$$
$$\pi(\varepsilon)\Big(\big(g^{-1}(\varepsilon)\big)^*q\wedge
g(\varepsilon)p_2\wedge\dots\wedge
g(\varepsilon)p_m\Big)=\big(g^{-1}(\varepsilon)\big)^*q\wedge
g(\varepsilon)p_2\wedge\dots\wedge g(\varepsilon)p_m$$
for all $q\in p^\perp$ and $p_1,\dots,p_m\in p$ since
$\big(g^{-1}(\varepsilon)\big)^*q\in\big(g(\varepsilon)p\big)^\perp$.
Taking derivatives, we obtain
$$\frac d{d\varepsilon}\Big|_{\varepsilon=0}
\pi'(\varepsilon)(p_1\wedge\dots\wedge
p_m)+\pi'\big[{\textstyle\bigwedge^m}p\big]
\sum\limits_{i=1}^mp_1\wedge\dots\wedge tp_i\wedge\dots\wedge
p_m=\sum\limits_{i=1}^mp_1\wedge\dots\wedge tp_i\wedge\dots\wedge p_m$$
and
$$\frac d{d\varepsilon}\Big|_{\varepsilon=0}\pi(\varepsilon)(q\wedge
p_2\wedge\dots\wedge p_m)+\pi\big[{\textstyle\bigwedge^mp}\big]\frac
d{d\varepsilon}\Big|_{\varepsilon=0}
\Big(\big(g^{-1}(\varepsilon)\big)^*q\wedge
g(\varepsilon)p_2\wedge\dots\wedge g(\varepsilon)p_m\Big)=$$
$$=\frac d{d\varepsilon}\Big|_{\varepsilon=0}
\Big(\big(g^{-1}(\varepsilon)\big)^*q\wedge
g(\varepsilon)p_2\wedge\dots\wedge g(\varepsilon)p_m\Big).$$
From $t^*q\in p$ and from
$$\frac d{d\varepsilon}\Big|_{\varepsilon=0}
\Big(\big(g^{-1}(\varepsilon)\big)^*q\wedge
g(\varepsilon)p_2\wedge\dots\wedge g(\varepsilon)p_m\Big)=-t^*q\wedge
p_2\wedge\dots\wedge p_m+\sum\limits_{i=2}^mq\wedge
p_2\wedge\dots\wedge tp_i\wedge\dots\wedge p_m,$$
we conclude that
$$\frac d{d\varepsilon}\Big|_{\varepsilon=0}
\pi'(\varepsilon)(p_1\wedge\dots\wedge
p_m)=\sum\limits_{i=1}^mp_1\wedge\dots\wedge tp_i\wedge\dots\wedge
p_m,\leqno{\bold{(2.7)}}$$
$$\frac d{d\varepsilon}\Big|_{\varepsilon=0}\pi(\varepsilon)(q\wedge
p_2\wedge\dots\wedge p_m)=-t^*q\wedge p_2\wedge\dots\wedge
p_m.\leqno{\bold{(2.8)}}$$

Let $X$ be a lifted field over a neighbourhood of $p$. Denote
$X(\varepsilon):=X\big(g(\varepsilon)p\big)$ and $s:=X(0)=X(p)$.
Define
$$E(\varepsilon):{\textstyle\bigwedge^m}V\to
{\textstyle\bigwedge^m}V,\qquad
v_1\wedge\dots\wedge v_m\mapsto\sum\limits_{i=1}^mv_1\wedge\dots\wedge
X(\varepsilon)v_i\wedge\dots\wedge v_m.$$
Clearly,
$E^mX(\varepsilon)=\pi(\varepsilon)E(\varepsilon)\pi'(\varepsilon)$. We
conclude from (2.7), (2.8), and $st=0$ that
$$\nabla_{E^mt}E^mX(p_1\wedge\dots\wedge p_m)=\Big(\frac
d{d\varepsilon}\Big|_{\varepsilon=0}\pi(\varepsilon)E(\varepsilon)
\pi'(\varepsilon)\Big)_{\bigwedge^mp}p_1\wedge\dots\wedge p_m=$$
$$=\pi(0)\Big(\frac
d{d\varepsilon}\Big|_{\varepsilon=0}\pi(\varepsilon)E(0)+\frac
d{d\varepsilon}\Big|_{\varepsilon=0}E(\varepsilon)+E(0)\frac
d{d\varepsilon}\Big|_{\varepsilon=0}\pi'(\varepsilon)\Big)
p_1\wedge\dots\wedge p_m=$$
$$=\pi(0)\Big(-\sum\limits_{i=1}^m p_1\wedge\dots\wedge t^*sp_i\wedge
\dots\wedge p_m+\sum\limits_{i=1}^mp_1\wedge\dots\wedge\frac
d{d\varepsilon}\Big|_{\varepsilon=0}X(\varepsilon)p_i\wedge\dots\wedge
p_m+$$
$$+\sum_{i\ne j}p_1\wedge\dots\wedge sp_i\wedge\dots\wedge
tp_j\wedge\dots\wedge p_m\Big)=
\sum\limits_{i=1}^mp_1\wedge\dots\wedge\pi[p]\frac
d{d\varepsilon}\Big|_{\varepsilon=0}X(\varepsilon)p_i\wedge\dots\wedge
p_m+B(s,t)p_1\wedge\dots\wedge p_m$$
(in the terms of the connection in
$\Gr_\Bbb K^0\big({k\choose m},\bigwedge^mV\big)$). In other words,
$$\nabla_{E^mt}E^mX=E^m\nabla_tX+B\big(X(p),t\big).$$
The first term is tangent to the image of the $m$-Pl\"ucker embedding
and the second one is orthogonal to~it $_\blacksquare$

\medskip

{\bf2.9.~Corollary.} {\sl The intrinsic connection is hermitian\/
{\rm(}pseudo-riemannian\/{\rm).}}

\medskip

{\bf Proof.} Taking $m=k$, the fact follows from Propositions 2.5, 2.6,
and [AGr, Proposition 4.3]
$_\blacksquare$

\medskip

{\bf2.10.~Corollary.} {\sl Let\/ $p\in\Gr_\Bbb K^0(k,V)$ and let\/
$t,t_1,t_2:p\to p^\perp$ be tangent vectors to\/ $\Gr_\Bbb K^0(k,V)$
at\/ $p$. The curvature tensor is given by}
$$R(t_1,t_2)t=tt_1^*t_2+t_2t_1^*t-tt_2^*t_1-t_1t_2^*t.$$

\medskip

{\bf Proof.} Since the above formula provides the curvature tensor in
the projective case [AGr, Subsection~4.4], it suffices to show that the
curvature tensors in $\Gr_\Bbb K^0(k,V)$ and in
$\Gr_\Bbb K^0\big({k\choose m},\bigwedge^mV\big)$ given by this formula
satisfy the Gauss equation (see [KoN, Proposition VII.4.1]) related to
the embedding $E^m$.

Let $t,t_1,t_2:p\to p^\perp$ be tangent vectors. Then, by Lemma 2.3,
$$E^mt(E^mt_1)^*E^mt_2(p_1\wedge\dots\wedge
p_m)=E^mt\sum\limits_{i=1}^mp_1\wedge\dots\wedge
t_1^*t_2p_i\wedge\dots\wedge p_m=$$
$$=\sum\limits_{i\ne j}p_1\wedge\dots\wedge tp_i\wedge\dots\wedge
t_1^*t_2p_j\wedge\dots\wedge
p_m+\sum\limits_{i=1}^mp_1\wedge\dots\wedge
tt_1^*t_2p_i\wedge\dots\wedge p_m.$$
for all $p_1,\dots,p_m\in p$. The last sum is exactly
$E^m(tt_1^*t_2)(p_1\wedge\dots\wedge p_m)$. Hence,
$$\big(E^mt(E^mt_1)^*E^mt_2-E^m(tt_1^*t_2)\big)(p_1\wedge\dots\wedge
p_m)=\sum_{i\ne j}p_1\wedge\dots\wedge tp_i\wedge\dots\wedge
t_1^*t_2p_j\wedge\dots\wedge p_m=B(t,t_1^*t_2).$$
Therefore, the Gauss equation takes the form\footnote{Strictly
speaking, we should take the (pseudo-)riemannian metric in the
equality. However, the Gauss equation turns out to be valid in a sense
which is even stronger than the hermitian one.}
$$\big\langle
E^mw,B(t,t_1^*t_2)+B(t_2,t_1^*t)-B(t,t_2^*t_1)-B(t_1,t_2^*t)\big\rangle
=\big\langle B(t_1,w),B(t_2,t)\big\rangle-\big\langle
B(t_2,w),B(t_1,t)\big\rangle,$$
where $w:p\to p^\perp$. So, it suffices to show that
$$(E^mw)^*B(t,t_1^*t_2)+(E^mw)^*B(t_2,t_1^*t)=
\big(B(t_1,w)\big)^*B(t_2,t),$$
$$(E^mw)^*B(t,t_2^*t_1)+(E^mw)^*B(t_1,t_2^*t)=
\big(B(t_2,w)\big)^*B(t_1,t).$$
We prove only the first identity. By Lemma 2.3,
$$(E^mw)^*B(t,t_1^*t_2)(p_1\wedge\dots\wedge p_m)=\sum_{i\ne
j}p_1\wedge\dots\wedge w^*tp_i\wedge\dots\wedge
t_1^*t_2p_j\wedge\dots\wedge p_m,$$
$$(E^mw)^*B(t_2,t_1^*t)(p_1\wedge\dots\wedge p_m)=\sum_{i\ne
j}p_1\wedge\dots\wedge w^*t_2p_i\wedge\dots\wedge
t_1^*tp_j\wedge\dots\wedge p_m,$$
and
$$(B(t_1,w))^*B(t_2,t)(p_1\wedge\dots\wedge
p_m)=\big(B(t_1,w)\big)^*\sum_{i\ne j}p_1\wedge\dots\wedge
t_2p_i\wedge\dots\wedge tp_j\wedge\dots\wedge p_m=$$
$$=\sum_{i\ne j}p_1\wedge\dots\wedge t_1^*t_2p_i\wedge\dots\wedge
w^*tp_j\wedge\dots\wedge p_m+\sum_{i\ne j}p_1\wedge\dots\wedge
w^*t_2p_i\wedge\dots\wedge t_1^*tp_j\wedge\dots\wedge p_m\
_\blacksquare$$

\medskip

{\bf2.11.~Corollary {\rm(compare to [BoN, Assertions 1--2])}.} {\sl
The\/ $m$-Pl\"ucker embedding is minimal.}

\medskip

{\bf Proof.} Let $e_1,\dots,e_k$ and $f_1,\dots,f_{n-k}$ be orthonormal
bases in $p$ and $p^\perp$. We define $t_{ij}e_j:=f_i$ and
$t_{ij}e_m:=0$ if $m\ne j$, getting in this way an orthonormal basis in
the tangent space at $p$. It is easy to see that $B(t_{ij},t_{ij})=0$.
It remains to apply [dCa, Definition 2.10]
$_\blacksquare$

\medskip

{\bf2.12.~Corollary {\rm(compare to [BoN, pp.~53 and 63])}.}
{\sl$\Gr_\Bbb K^0(k,V)$ is Einstein. The corresponding constant is\/
$n-2$ in the case of\/ $\Bbb K=\Bbb R$ and\/ $2n$ in the case of\/
$\Bbb K=\Bbb C$, where\/ $n=\dim_\Bbb KV$.}

\medskip

{\bf Proof.} We use the following elementary fact: Let $T:V\to V$ be an
$\Bbb R$-linear map. Then $\tr_\Bbb RT=2\Re\tr_\Bbb CT$ if $T$ is
$\Bbb C$-linear and $\tr_\Bbb RT=0$ if $T$ is $\Bbb C$-antilinear.

The Ricci tensor is given by
$\ricci(t_1,t):=\tr\big(t_2\mapsto R(t_1,t_2)t\big)$, where
$t,t_1,t_2:p\to p^\perp$. Considering each term of the curvature tensor
in Corollary 2.10, it is easy to see that
$$\tr(t_2\mapsto tt_1^*t_2)=k\tr(tt_1^*)=k\tr(t^*t_1),\quad
\tr(t_2\mapsto t_2t_1^*t)=(n-k)\tr(t_1^*t)=(n-k)\tr(t^*t_1),$$
$$\tr(t_2\mapsto tt_2^*t_1)=\tr(t_2\mapsto t_1t_2^*t)=\tr(t^*t_1)$$
in the case of $\Bbb K=\Bbb R$, and that
$${\tr}_\Bbb C(t_2\mapsto tt_1^*t_2)=k\tr(tt_1^*),\quad
{\tr}_\Bbb C(t_2\mapsto t_2t_1^*t)=(n-k)\tr(t_1^*t),$$
$${\tr}_\Bbb R(t_2\mapsto tt_1^*t_2)=2k\Re\tr(t^*t_1),\quad
{\tr}_\Bbb R(t_2\mapsto t_2t_1^*t)=2(n-k)\Re\tr(t^*t_1),$$
$${\tr}_\Bbb R(t_2\mapsto tt_2^*t_1)={\tr}_\Bbb R(t_2\mapsto
t_1t_2^*t)=0$$
in the case of $\Bbb K=\Bbb C$
$_\blacksquare$

\medskip

{\bf2.13.~Generic geodesics.} Let $p\in\Gr_\Bbb K^0(k,V)$ and let
$t\in\Lin_\Bbb K(p,p^\perp)\subset\Lin_\Bbb K(V,V)$ be a tangent vector
at $p$. We are going to describe the geodesic determined by $t$ in the
generic case, i.e., when there exists an orthonormal basis
$p_1,\dots,p_k$ in $p$ formed by nonisotropic eigenvectors of the
self-adjoint map $t^*t:p\to p$ (if $\Bbb K=\Bbb C$, this means that
$t^*t:p\to p$ has no isotropic eigenvectors).

The eigenvalues $\lambda_1,\dots,\lambda_k$ corresponding to
$p_1,\dots,p_k$ are real. Put $W_j:=\Bbb Rp_j+\Bbb Rtp_j$. The $W_j$'s
are pairwise orthogonal because the $tp_j$'s are pairwise orthogonal.
Being restricted to $W_j$, the form is real and does not vanish. So,
$W_j$ provides a geodesic $\G_j\subset\Bbb P_\Bbb KV$ if $tp_j\ne0$. By
[AGr, Lemma 2.1], $\G_j$~is respectively spherical, hyperbolic, or
euclidean exactly when $\lambda_j>0$, $\lambda_j<0$, or $\lambda_j=0$.
(If $tp_j=0$, $\G_j$~is a single point in $\Bbb P_\Bbb KV$.)

Let $t_j$ be the tangent vector to $\G_j$ at $p_j$ given by
$t_j:p_j\mapsto tp_j$. Every geodesic $\G_j$ admits a local {\it
uniformly\/} parameterized lift $p_j(s)$ to $V$ with respect to $t_j$.
This means that the tangent vector $p_j(s)\mapsto\dot p_j(s)$ at
$p_j(s)$ is the parallel displacement of $t_j$ from $p_j(0)=p_j$ to
$p_j(s)$ along $\G_j$ (in particular,
$\dot p_j(s)\in p_j(s)^\perp\cap W_j$) and that
$\big\langle p_j(s),p_j(s)\big\rangle$ is constant in $s$. If $\G_j$ is
noneuclidean, such a parameterization is readily obtainable from those
in [AGr, Subsection 3.2]. In the euclidean case, $p_j(s):=p_j+stp_j$ is
the desired parameterization [AGr, Corollary 5.9]. Note that
$\ddot p_j(s)\in\Bbb Rp_j(s)$. This is obvious in the euclidean case
and is otherwise implied by the fact that
$\big\langle\dot p_j(s),\dot p_j(s)\big\rangle$ is constant and
$\ddot p_j(s)\in W_j$.

As in [AGoG, Section 2], we fix a $k$-dimensional $\Bbb K$-vector space
$P$. Let $b_1,\dots,b_k\in P$ be a basis and let $p(s):P\to V$ be the
linear map given by the rule $p(s):b_j\mapsto p_j(s)$.

\medskip

{\bf2.14.~Lemma.} {\it The curve\/ $\G:s\mapsto p(s)$ is a geodesic
in\/ $\Gr_\Bbb K^0(k,V)$ and\/ $t$ is its  tangent vector at\/ $p$.}

\medskip

{\bf Proof.} The tangent vector to $\G$ at $p(s)$ is given by the
linear map
$t(s)\in\Lin_\Bbb K\big(p(s),p(s)^\perp\big)\subset\Lin_\Bbb K(V,V)$,
$t(s):p_j(s)\mapsto\dot p_j(s)$, because $\dot p_j(s)\in p(s)^\perp$
and the $W_j$'s are pairwise orthogonal.

In the definition of $\nabla$, taking the derivative of
$X\big(c(\varepsilon)\big)$ at $\varepsilon=0$, where
$c(\varepsilon):=(1+\varepsilon t)p$, amounts to taking the derivative
of $X\big(p(s)\big)$ at $s$ because $\dot c(0)=\dot p(s)$. Therefore,
$\nabla_{\dot{\G}(s)}\dot{\G}(s)=\pi\big[p(s)\big]\dot
t(s)\pi'\big[p(s)\big]$. Taking the derivative of
$t(s)p_j(s)=\dot p_j(s)$, we obtain
$\dot t(s)p_j(s)+t(s)\dot p_j(s)=\ddot p_j(s)$. Since
$t(s)\big(p(s)^\perp\big)=0$ and $\dot p_j(s)\in p(s)^\perp$, we have
$\pi\big[p(s)\big]\dot t(s)p_j(s)=\pi\big[p(s)\big]\ddot p_j(s)=0$ due
to $\ddot p_j(s)\in\Bbb Rp_j(s)$
$_\blacksquare$

\medskip

We call $\G_j$ a {\it spine\/} of $\G$. We may interpret a point
$\G(s)$ as a linear subspace in $\Bbb P_\Bbb KV$ spanned by the
$p_j(s)'s$. Moving along the geodesic $\G$ in $\Gr_\Bbb K^0(k,V)$ is
the same as moving along the spines with velocities given
by\footnote{Well, when an euclidean spine is involved the situation is
more subtle.}
$\sqrt{|\lambda_j|}$. The equality $tp_j=0$ says that $\G_j$ is a point
fixed during the movement.

A generic tangent vector $t$ provides a choice of a basis formed by the
eigenvectors of $t^*t$. In other words, if $2k\le n$, the intention of
moving in some generic direction automatically chooses a certain
reference frame.

\medskip

{\bf2.15.~Comments and questions.} Many of the above facts admit a form
not involving the hermitian metric.

\smallskip

$\bullet$ The first formula displayed in the proof of Proposition 2.5
says that $(E^mt_1)^*E^mt_2=E^m(t_1^*t_2)$.

\smallskip

$\bullet$ The Gauss equation in Corollary 2.10 follows from the much
simpler one
$(E^mw)^*B(t,t_1^*t_2)+(E^mw)^*B(t_2,t_1^*t)=\big(B(t_1,w)\big)^*
B(t_2,t)$.

\smallskip

$\bullet$ The proof of minimality actually does not require the
self-adjoint operator $S_\eta$ from [dCa, Definition~2.10].

\smallskip

$\bullet$ What is the geometrical meaning of the other two
symmetrizations of the trilinear product $tt_2^*t_1$ ?

\smallskip

$\bullet$ What about other functors in place of $\bigwedge^m$ ?

\bigskip

\centerline{\bf3.~Convexity of some real hyperbolic polyhedra}

\medskip

This section illustrates how grassmannians appear in a typical
situation that does not seem to involve them at the first glance. Here
we deal with the real hyperbolic geometry $\Bbb H_\Bbb R^4$, that is,
with $\Bbb P_\Bbb RV$, where $V$ is an $\Bbb R$-vector space and the
form has signature $++++-$. (The calculus in what follows may seem a
little bit concise. On the other hand, it requires no specific
knowledge in the area.)

A known problem on real hyperbolic disc bundles is to find the greatest
value of $|e/\chi|$, where $e$ stands for the Euler number of the
bundle and $\chi$, for the Euler characteristic of the base closed
surface [GLT]. By now, the best value $|e/\chi|=1/2$ [Kui], [Luo] is
obtained via constructing a fundamental polyhedron without faces of
codimension $>2$ that is strongly convex in the sense that its disjoint
faces lie in disjoint totally geodesic hypersurfaces. It is worthwhile
trying polyhedra that are convex in the usual sense.

Such a polyhedron can be described in the terms of a finite number of
positive points $p_1,\dots,p_n\in\Bbb P_\Bbb RV$. The face $F_i$ is
a {\it segment\/} in the hyperplane $H_i:=p_i^\perp\cap\overline\B V$,
i.e., the part of $H_i$ between the disjoint planes $E_{i-1}$ and
$E_i$, where
$E_i:=F_i\cap F_{i+1}=\Span(p_i,p_{i+1})^\perp\cap\overline\B V$ for
all $i$ (the indices are modulo~$n$). In the terms of the Gram matrix
$U(p_1,\dots,p_n):=[u_{ij}]$, $u_{ij}:=\langle p_i,p_j\rangle$,
assuming that $u_{ii}=1$, the~strong convexity means
$|u_{i(i+1)}|<1<|u_{ij}|$ for all $j\ne i-1,i,i+1$. In what follows, we
obtain a criterion for the usual convexity.

It is convenient to use the following notation:
$$\langle i_1i_2,j_1j_2\rangle:=\det\pmatrix
u_{i_1j_1}&u_{i_1j_2}\\u_{i_2j_1}&u_{i_2j_2}\endpmatrix,\qquad\langle
i_1i_2i_3,j_1j_2j_3\rangle:=\det\pmatrix
u_{i_1j_1}&u_{i_1j_2}&u_{i_1j_3}\\u_{i_2j_1}&u_{i_2j_2}&u_{i_2j_3}\\
u_{i_3j_1}&u_{i_3j_2}&u_{i_3j_3}\endpmatrix.\leqno{\bold{(3.1)}}$$

The fact that $H_i\cap H_{i+1}\ne\varnothing$ can be written as
$\big\langle i(i+1),i(i+1)\big\rangle>0$.
 The fact that $E_{i-1}$ and
$E_i$ are disjoint is equivalent to
$\Span(p_{i-1},p_i,p_{i+1})^\perp\cap\overline\B V=\varnothing$, i.e.,
to $\big\langle(i-1)i(i+1),(i-1)i(i+1)\big\rangle<0$ by Sylvester's
criterion.

\medskip

{\bf3.2.~Lemma.} {\sl The segment $F_i$ can be described as}
$$F_i=\Big\{x\in H_i\mid\big\langle(i-1)i,i(i+1)\big\rangle\langle
x,p_{i-1}\rangle\langle p_{i+1},x\rangle\ge0\Big\}.$$

{\bf Proof}. During the proof, we deal only with the points
$p_{i-1},p_i,p_{i+1}$. We change these points keeping $E_{i-1},F_i,E_i$
the same. The expression $\big\langle(i-1)i,i(i+1)\big\rangle$ does not
change if we substitute $p_{i-1}$ and $p_{i+1}$ respectively by
$p_{i-1}+r_1p_i$ and $p_{i+1}+r_2p_i$, $r_1,r_2\in\Bbb R$. Also,
$\big\langle(i-1)i,i(i+1)\big\rangle\langle x,p_{i-1}\rangle\langle
p_{i+1},x\rangle$
does not change if we alter the sign of $p_{i-1}$. So, we can assume
that
$u_{(i-1)i}=u_{i(i+1)}=0$, $u_{(i-1)(i-1)}=u_{ii}=u_{(i+1)(i+1)}=1$,
and $u_{(i-1)(i+1)}\ge0$. It follows from
$\big\langle(i-1)i(i+1),(i-1)i(i+1)\big\rangle<0$ that
$u_{(i-1)(i+1)}>1$. The closed $3$-ball $H_i$ is fibred over the
hyperbolic geodesic $\G_i:=\Span(p_{i-1},p_{i+1})$ by the closed discs
$S_p:=\Span(p,p_i)^\perp\cap\overline\B V$ called {\it slices,}
$p\in\G_i\setminus\overline\B V$. The end slices $E_{i-1}$ and $E_i$ of
$F_i$ correspond to $p=p_{i-1}$ and $p=p_{i+1}$. Since
$u_{(i-1)(i+1)}>0$, the segment $F_i$ is formed by the slices $S_p$
with $p=(1-t)p_{i-1}+tp_{i+1}$, $t\in[0,1]$. Note that
$\Span(p_{i-1},p_{i+1})=\Span({^{p_{i+1}}}p_{i-1},{^{p_{i-1}}}p_{i+1})$
because $u_{(i-1)(i+1)}>1$.

Let $x\in H_i$. Then
$x=w-t_1{^{p_{i+1}}}p_{i-1}+t_2{^{p_{i-1}}}p_{i+1}$ for suitable
$w\in\Span(p_{i-1},p_{i+1})^\perp$, $t_1,t_2\in\Bbb R$, $t_1\ge0$. We
have
$$\big\langle(i-1)i,i(i+1)\big\rangle\langle x,p_{i-1}\rangle\langle
p_{i+1},x\rangle=u_{(i-1)(i+1)}(u_{(i-1)(i+1)}^2-1)^2t_1t_2$$
and $\langle t_2p_{i-1}+t_1p_{i+1},x\rangle=0$. It follows from
$x\in\overline\B V$ that $t_2p_{i-1}+t_1p_{i+1}\notin\overline\B V$ and
that $t_1\ne0$ or $t_2\ne0$. So, $x\in S_{t_2p_0+t_0p_2}$ and the claim
easily follows
$_\blacksquare$

\medskip

In the sequel, we frequently use the above decomposition of $H_i$ into
slices over the hyperbolic geodesic~$G_i$.

\smallskip

The usual convexity is equivalent to the condition
$F_i\cap H_j=\varnothing$ for $j\ne i-1,i,i+1$. We fix $i$ and $j$ and
express this condition by considering the following cases:

\smallskip

$\bullet$ $\langle ij,ij\rangle<0$. This implies
$H_i\cap H_j=\varnothing$, hence, $F_i\cap H_j=\varnothing$.

\smallskip

$\bullet$ $\langle ij,ij\rangle=0$. First, we require $p_j\ne p_i$
(implied by $F_i\cap H_j=\varnothing$). Under these conditions,
the~isotropic point $u_{ii}p_j-u_{ji}p_i$ is the only point in
$\Span(p_i,p_j)^\perp\cap\overline\B V$. By Lemma 3.2, the condition
$F_i\cap H_j=\varnothing$ is equivalent to
$$\big\langle ij,(i-1)i\big\rangle\big\langle (i-1)i,i(i+1)\big
\rangle\big\langle i(i+1),ij\big\rangle>0.\leqno{\bold{(3.3)}}$$
It obviously implies that $p_j\ne p_i$.

\smallskip

$\bullet$ $\langle ij,ij\rangle>0$. Define
$$q_1:=\frac{u_{ii}p_{i-1}-u_{(i-1)i}p_i}
{\sqrt{u_{ii}\big\langle(i-1)i,(i-1)i\big\rangle}},\quad
q_2:=\frac{u_{ii}p_{i+1}-u_{(i+1)i}p_i}
{\sqrt{u_{ii}\big\langle(i+1)i,(i+1)i\big\rangle}},\quad
q_3:=\frac{u_{ii}p_j-u_{ji}p_i}
{\sqrt{u_{ii}\langle ij,ij\rangle}},$$
and $v_{kl}:=\langle q_k,q_l\rangle$. It is easy to see that
$q_k\in p_i^\perp$ and $v_{kk}=1$ for all $k$. The facts that
$\Span(q_1,q_2,p_i)=\Span(p_{i-1},p_i,p_{i+1})$ has signature $++-$ and
that $p_i$ is positive imply $|v_{12}|>1$. The slices of $F_i$ have the
form $S_{q(t)}$, where
$$q(t):=(1-t)q_1+\sigma tq_2,\qquad t\in[0,1],$$
and $\sigma:=\displaystyle\frac{v_{12}}{|v_{12}|}$. The condition
$F_i\cap H_j=\varnothing$ is equivalent to the requirement that
$\Span\big(q(t),q_3\big)$ has signature $+-$ for all $t\in[0,1]$. It
can be written as
$$f(t):=t^2\big((v_{13}-\sigma v_{23})^2+2|v_{12}|-2\big)
-2t\big(v_{13}^2-\sigma v_{13}v_{23}+|v_{12}|-1\big)+v_{13}^2-1>0$$
by Sylvester's criterion.

Writing $f(t)=t^2a-2tb+c$, we have $a>0$, $f(0)=c=v_{13}^2-1$, and
$f(1)=v_{23}^2-1$. The polynomial $f(t)$ attains its minimum at
$t=b/a$. Clearly, $f(b/a)>0$ if and only if $ac>b^2$. Hence, the
condition $F_i\cap H_j=\varnothing$ is equivalent to
$v_{13}^2,v_{23}^2>1$ and $0<b<a\Longrightarrow ac>b^2$. One readily
verifies that
$$ac-b^2=1+2v_{12}v_{23}v_{31}-v_{12}^2-v_{23}^2-v_{31}^2=\det
U(q_1,q_2,q_3)$$
and that $0<b<a$ is equivalent to
$$v_{13}^2-1>\sigma(v_{13}v_{32}-v_{12}),\qquad
v_{23}^2-1>\sigma(v_{13}v_{32}-v_{12}).$$
The inequality $ac>b^2$ is impossible because $\Span(q_1,q_2,q_3)$
contains a negative point belonging to $\G_i=\Span(q_1,q_2)$.
Therefore, $F_i\cap H_j=\varnothing$ is equivalent to
$v_{13}^2,v_{23}^2>1$ and $v_{13}^2-1\le\sigma(v_{13}v_{32}-v_{12})$ or
$v_{23}^2-1\le\sigma(v_{13}v_{32}-v_{12})$. Either of the last two
inequalities implies $\sigma(v_{13}v_{32}-v_{12})>0$, that is,
$\sigma v_{13}v_{32}>|v_{12}|$, i.e.,
$$v_{12}v_{23}v_{31}>v_{12}^2.\leqno{\bold{(3.4)}}$$
Clearly, (3.4) implies $\sigma(v_{13}v_{32}-v_{12})>0$. Assuming that
(3.4) is true, we can rewrite the condition
$v_{13}^2-1\le\sigma(v_{13}v_{32}-v_{12})$ or
$v_{23}^2-1\le\sigma(v_{13}v_{32}-v_{12})$ in the form
$$(v_{13}^2-1)^2\le(v_{13}v_{32}-v_{12})^2\qquad\text{or}\qquad
(v_{23}^2-1)^2\le(v_{13}v_{32}-v_{12})^2.\leqno{\bold{(3.5)}}$$
In fact, the meaning of the inequalities $v_{13}^2,v_{23}^2>1$ is that
$\Span(p_{i-1},p_i,p_j)$ and $\Span(p_{i+1},p_i,p_j)$ have signature
$++-$, that is,
$$\big\langle(i-1)ij,(i-1)ij\big\rangle<0,\qquad\big\langle
i(i+1)j,i(i+1)j\big\rangle<0.$$
Under these conditions, (3.5) takes the form
$$\min(v_{13}^2-1,v_{23}^2-1)\le|v_{13}v_{32}-v_{12}|
\leqno{\bold{(3.6)}}$$
By straightforward calculus, we have
$$v_{12}=-\frac{\big\langle(i-1)i,i(i+1)\big\rangle}{\sqrt{\big\langle
(i-1)i,(i-1)i\big\rangle\big\langle i(i+1),i(i+1)\big\rangle}},\qquad
v_{23}=\frac{\big\langle i(i+1),ij\big\rangle}{\sqrt{\big\langle
i(i+1),i(i+1)\big\rangle\langle ij,ij\rangle}},$$
$$v_{13}=-\frac{\big\langle(i-1)i,ij\big\rangle}
{\sqrt{\big\langle(i-1)i,(i-1)i\big\rangle\langle ij,ij\rangle}}.$$
Hence, (3.4) takes the form
$$\frac{\big\langle(i-1)i,ij\big\rangle\big\langle
ij,i(i+1)\big\rangle}{\big\langle(i-1)i,i(i+1)\big\rangle}>\langle
ij,ij\rangle.\leqno{\bold{(3.7)}}$$
Note that (3.7) is equivalent to (3.3) in the case of
$\langle ij,ij\rangle=0$ because
$\big\langle(i-1)i,i(i+1)\big\rangle=0$ would imply $v_{12}=0$, that
is, $\langle{^{p_i}}p_{i-1},{^{p_i}}p_{i+1}\rangle=0$, contradicting
$E_{i-1}\cap E_i=\varnothing$.

Since
$$v_{13}v_{32}-v_{12}=\frac{\big\langle(i-1)i,i(i+1)\big\rangle\langle
ij,ij\rangle-\big\langle(i-1)i,ij\big\rangle\big\langle
i(i+1),ij\big\rangle}{\langle
ij,ij\rangle\sqrt{\big\langle(i-1)i,(i-1)i\big\rangle\big\langle
i(i+1),i(i+1)\big\rangle}}=$$
$$=\frac{u_{ii}\big\langle(i-1)ij,i(i+1)j\big\rangle}{\langle
ij,ij\rangle\sqrt{\big\langle(i-1)i,(i-1)i\big\rangle\big\langle
i(i+1),i(i+1)\big\rangle}},$$
$$v_{13}^2-1=\frac{\big\langle(i-1)i,ij\big\rangle\big\langle(i-1)i,
ij\big\rangle-\big\langle(i-1)i,(i-1)i\big\rangle\langle
ij,ij\rangle}{\big\langle(i-1)i,(i-1)i\big\rangle\langle
ij,ij\rangle}=-\frac{u_{ii}\big\langle(i-1)ij,(i-1)ij\big\rangle}
{\langle ij,ij\rangle\big\langle(i-1)i,(i-1)i\big\rangle},$$
$$v_{23}^2-1=\frac{\big\langle i(i+1),ij\big\rangle\big\langle
i(i+1),ij\big\rangle-\big\langle i(i+1),i(i+1)\big\rangle\langle
ij,ij\rangle}{\big\langle i(i+1),i(i+1)\big\rangle\langle
ij,ij\rangle}=-\frac{u_{ii}\big\langle
i(i+1)j,i(i+1)j\big\rangle}{\langle ij,ij\rangle\big\langle
i(i+1),i(i+1)\big\rangle},$$
$u_{ii}>0$, and $\langle ij,ij\rangle>0$, (3.6) takes the form
$$\frac{\Big|\big\langle(i-1)ij,i(i+1)j\big\rangle\Big|}
{\sqrt{\big\langle(i-1)i,(i-1)i\big\rangle\big\langle
i(i+1),i(i+1)\big\rangle}}+\max\bigg(\frac{\big\langle(i-1)ij,(i-1)
ij\big\rangle}{\big\langle(i-1)i,(i-1)i\big\rangle},\frac{\big\langle
i(i+1)j,i(i+1)j\big\rangle}{\big\langle
i(i+1),i(i+1)\big\rangle}\bigg)\ge0.\leqno{\bold{(3.8)}}$$

It follows from $E_{i-1}\subset F_i$ and $F_i\cap H_j=\varnothing$
that $E_{i-1}\cap H_j=\varnothing$ for $j\ne i-1,i,i+1$. In other
words, $H_{i-1}\cap H_i\cap H_j=\varnothing$, that is,
$\big\langle(i-1)ij,(i-1)ij\big\rangle<0$.

\smallskip

Summarizing, we arrive at the

\medskip

{\bf3.9.~Criterion of convexity.} {\sl The polyhedron formed by
segments of hyperplanes given by\/ $p_1,\dots,p_n\allowmathbreak\in V$
is convex\/ {\rm(}hence, simple\/{\rm)} if and only if the following
conditions written in the terms of\/ {\rm(3.1)}, where\/
$u_{ij}:=\langle p_i,p_j\rangle$, hold\/ {\rm(}the indices are modulo\/
$n${\rm):}

\smallskip

$\bullet$ The inequalities\/ $u_{ii}>0$ are valid for all\/ $i$.

$\bullet$ The inequalities\/ $\big\langle(i-1)i,(i-1)i\big\rangle>0$
and\/ $\big\langle(i-1)ij,(i-1)ij\big\rangle<0$ are valid for all\/
$j\ne i-1,i$.

$\bullet$ The inequalities\/ {\rm(3.7)} are valid for all\/
$j\ne i-1,i,i+1$ such that\/ $\langle ij,ij\rangle\ge0$.

$\bullet$ The inequalities\/ {\rm(3.8)} are valid for all\/
$j\ne i-1,i,i+1$ such that\/ $\langle ij,ij\rangle>0$
$_\blacksquare$}

\medskip

Note that
$$\langle i_1i_2,j_1j_2\rangle=\langle
g_{i_1i_2},g_{j_1j_2}\rangle,\qquad\langle
i_1i_2i_3,j_1j_2j_3\rangle=\langle
g_{i_1i_2i_3},g_{j_1j_2j_3}\rangle,$$
where $g_{i_1i_2}:=p_{i_1}\wedge p_{i_2}$ and
$g_{i_1i_2i_3}:=p_{i_1}\wedge p_{i_2}\wedge p_{i_3}$ represent
respectively
$\bigwedge^2\Span(p_{i_1},p_{i_2})\in\Bbb P_\Bbb R\bigwedge^2V$ and
$\bigwedge^3\Span(p_{i_1},p_{i_2},p_{i_3})\in\Bbb P_\Bbb
R\bigwedge^3V$.
So, Criterion 3.9 deals with the usual projective invariants.

\bigskip

\centerline{\bf4.~References}

\medskip

[AGoG] S.~Anan$'$in, E.~C.~B.~Gon\c calves, C.~H.~Grossi, {\it
Grassmannians and conformal structure on absolutes,} preprint
http://arxiv.org/abs/0907.4469

[AGr] S.~Anan$'$in, C.~H.~Grossi, {\it Coordinate-free classic
geometries,} preprint\newline
http://arxiv.org/abs/math/0702714

[BoN] A.~A.~Borisenko, Yu.~A.~Nikolaevskii, {\it Grassmann manifolds
and the Grassmann image of submanifolds,} Russian Math. Surveys,
No.~46 (1991), 45--94

[dCa] M.~P.~do Carmo, {\it Riemannian Geometry,} Mathematics: Theory
\&{} Applications, Birkh\"auser, Boston, 1992. vi+300 pp.

[GLT] M.~Gromov, H.~B.~Lawson Jr., W.~Thurston, {\it Hyperbolic\/
$4$-manifolds and conformally flat\/ $3$-manifolds,} Inst.~Hautes
\'Etudes Sci.~Publ.~Math., No.~68 (1988), 27--45

[GuK] B.~Guilfoyle, W.~Klingenberg, {\it Proof of the Caratheodory
Conjecture by Mean Curvature Flow in the Space of Oriented Affine
Lines,} preprint http://arxiv.org/abs/0808.0851

[KNo] S.~Kobayashi, K.~Nomizu, {\it Foundations of Differential
Geometry,} Volume II, Interscience Tracts in Pure and Applied
Mathematics. Interscience Publishers. John Wiley \&{} Sons, 1969.
xv+470 pp.

[Kui] N.~H.~Kuiper, {\it Hyperbolic\/ $4$-manifolds and tessellations,}
Inst.~Hautes \'Etudes Sci.~Publ.~Math., No.~68 (1988), 47--76

[Luo] F.~Luo, {\it Constructing conformally flat structures on some
Seifert fibred\/ $3$-manifolds,} Math.~Ann. {\bf294} (1992), No.~3,
449--458

\enddocument